\documentclass[english,prl, twocolumn]{revtex4-1}

\usepackage{amsmath}
\usepackage{amssymb}
\usepackage{graphicx}

\usepackage{bm}

\usepackage[colorlinks=true,
            linkcolor=blue,
            urlcolor=blue,
            citecolor=blue]{hyperref}

\renewcommand{\vec}[1]{\bm{#1}}
\newcommand{\mat}[1]{\bm{#1}}
\newcommand{\map}{\vec{F}}

\newcommand{\koopman}{\mathcal{K}}
\newcommand{\state}{{\vec{x}}}
\newcommand{\manifold}{\mathcal{M}}

\newcommand{\dictionary}{\mat{\psi}}

\newcommand{\edmdK}{\mat{{K}}}
\newcommand{\obsa}{\vec{\tilde x}}
\newcommand{\obsaf}{\vec{\tilde y}}
\newcommand{\Ma}{\tilde{M}}
\newcommand{\obsfuna}{\vec{\tilde g}}

\newcommand{\obsb}{\vec{\hat x}}
\newcommand{\obsbf}{\vec{\hat y}}
\newcommand{\obsfunb}{\vec{\hat g}}
\newcommand{\Mb}{\hat{M}}

\begin{document}

\title{Data Fusion via Intrinsic Dynamic Variables:\\
An Application of Data-Driven Koopman Spectral Analysis}

\author{Matthew O. Williams}

\email{mow2@princeton.edu}
\affiliation{Program in Applied and Computational Mathematics (PACM), Princeton University, NJ 08544, USA.}

\author{Clarence W. Rowley}
\affiliation{Department of Mechanical and Aerospace Engineering, Princeton University, NJ 08544, USA.}

\author{Igor Mezi\'c}
\affiliation{Department of Mechanical Engineering, University of California, Santa Barbara, CA 93106, USA.}

\author{Ioannis G. Kevrekidis}
\affiliation{Department of Chemical and Biological Engineering \& PACM, Princeton University, NJ 08544, USA.}

\begin{abstract}
We demonstrate that numerically computed approximations of Koopman
eigenfunctions and eigenvalues create a natural framework for data fusion in
applications governed by nonlinear evolution laws.
This is possible because the eigenvalues of the Koopman operator are
invariant to invertible transformations of the system state, so that the values
of the Koopman eigenfunctions serve as a set of \emph{intrinsic coordinates}
that can be used to map between different observations (e.g., measurements
obtained through different sets of sensors) of the same fundamental behavior.
The measurements we wish to merge can also be nonlinear, but must be ``rich
enough'' to allow (an effective approximation of) the state to be reconstructed
from a single set of measurements.
This approach requires independently obtained sets of data that capture the
evolution of the heterogeneous measurements and a single pair
of ``joint'' measurements taken at one instance in time.
Computational approximations of eigenfunctions and their corresponding
eigenvalues from data are accomplished using Extended Dynamic Mode
Decomposition.
We illustrate this approach on measurements of spatio-temporal oscillations of
the FitzHugh-Nagumo PDE, and show how to fuse point measurements with principal
component measurements, after which either set of measurements can be used to estimate the other set.
\end{abstract}
\maketitle

In many applications, our understanding of a system comes from sets of partial measurements (functions
of the system state) rather than observations of the full state.
Linking these  heterogeneous partial measurements (different sensors, different measurement times)
is  the objective of a broad collection of techniques referred to as data fusion methods~\cite{Hall1997,Pohl1998,Yocky1995}.
 Since the system state can itself be a measurement, these methods also
 encompass traditional techniques for state estimation such as Kalman
 filters~\cite{Simon2006,Stengel2012,bishop2001introduction}  and stochastic estimation methods~\cite{Adrian1988,Adrian1994,Guezennec1989,Ho1964,Bonnet1994,Druault2010,Murray2007,Naguib2001}.
 One might subdivide these approaches into ``dynamic'' methods, which require models of the underlying
 evolution~\cite{Simon2006,Stengel2012,bishop2001introduction}, and ``static'' methods~\cite{Fieguth2010},
 which do not require dynamical information but, in general, need more extensive
 measurements to be successful.
Other, more recent, methods such as the gappy Proper Orthogonal Decomposition~\cite{Willcox2006},
nonlinear intrinsic variables~\cite{Dsilva2013}, or compressed sensing-based methods~\cite{Bright2013}
can be thought of as solving, in effect, the same problem, but make different
assumptions about the nature of the underlying system and the type of measurements used.
Though the exact details and assumptions differ, the overarching goal in data
fusion is to develop a mapping from one type of measurements to another type.

In this manuscript, we propose a method that generates such a mapping with the help of a set of
{\it intrinsic coordinates}; these coordinates are based on the (computationally approximated)
eigenfunctions of the Koopman operator~\cite{Koopman1931,Koopman1932,Mauroy2013,Mauroy2012,Budivsic2012}.
To merge measurements from heterogeneous sources,
the algorithm requires data in the form of time series from each source, and a
set of ``joint'' measurements (i.e., measurements known to correspond to the same underlying state).
Each of these individual measurements must be ``rich enough'' so that
the system state (or a quantity effectively equivalent to the state such as the
leading Principal Components~\cite{Lee2007}) can be recovered using a single
set of measurements, but this limitation can likely be overcome by enriching
measurement sets that are not ``rich enough'' through the use of time
delays~\cite{Chorin2000,Chorin2002,juang1994applied}.
The benefits of this approach are that it is naturally applicable to
nonlinear systems and sensors, and minimizes the number of joint measurements
required; in many systems, only a single pair is needed.

{\bf Problem formulation.} 
Suppose we have {\em two} different sets of measurements, generated by two
different sets of (heterogeneous) sensors observing the same fundamental
behavior (the evolution of the state $\vec{x}(t)$).
Let $\obsa=\obsfuna(\vec{x})$ denote a measurement of this state obtained with the first set of sensors, and $\obsb=\obsfunb(\vec{x})$  a measurement with the second set; ``a measurement" is, in general, a vector valued observable obtained at a single moment in time by one of our two collections of sensors.
Here, $\obsfuna:\manifold\to\tilde{\manifold}$ and
$\obsfunb:\manifold\to\hat{\manifold}$ are the functions that map from the
instantaneous system state ($\vec x \in\manifold$) to the corresponding instantaneous measurements.
We record time series of such measurements from each of our sets of sensors, and
divide each of the time series into sets of measurement  pairs,
$\{(\obsa_m,\obsaf_m)\}_{m=1,\ldots,\Ma}$ and $\{(\obsb_{m},\obsbf_{m})\}_{m=1,\ldots,\Mb}$, where $\obsa_{m}$ (resp. $\obsb_{m}$) is the $m$-th measurement, and $\obsaf_{m}$ (resp. $\obsbf_{m}$) is its value after a single sampling interval.
These time-series can be obtained independently, and the total number of measurements, $\Ma$ and $\Mb$, can differ.
For simplicity, we assume the sampling interval, $\Delta t$, is the same for both data sets; this too can vary with only a slight modification to the algorithm.
The only requirements we place on these data sets are that: (i) the
data they contain are generated by the same system, (ii) the state can, in
principle, be determined using a snapshot from either collection of
measurements, and (iii) that the sampling interval remains constant {\em within} a data set.
The (required) {\em joint data set} is denoted as  $\left\{ (\obsa_{m},\obsb_{m})\right\}_{m=1}^{M_\text{joint}} $;
the subscripts denote the $m$-th measurement in the joint data set, and $M_\text{joint}=1$ is the number of data pairs.
This approach is applicable to an arbitrary number of different measurements;
the restriction to two is solely for simplicity of presentation.

{\bf The Koopman operator.}
The crux of our approach is the use of these time-series data to computationally approximate the leading eigenfunctions of the Koopman operator~\cite{Koopman1931,Koopman1932,Budivsic2012,Mauroy2013,Mauroy2013a,Mezi2013},
thus generating a mapping from measurement space to {\it an intrinsic variable space}.
The Koopman operator is defined for a specific autonomous dynamical system, whose evolution law we denote as
${\vec{x}\mapsto\vec{F}(\vec{x}),}$
where $\vec{x}\in\manifold\subseteq\mathbb{R}^{N}$ is the system state, $\vec{F}:\manifold\to\manifold$ is the evolution operator, and $n\in\mathbb{N}$ is discrete time.
The action of the Koopman operator is
\begin{equation}
(\koopman\psi)(\state)=(\psi\circ\map)(\state)=\psi(\map(\state)),\label{eq:koopman}
\end{equation}
where $\psi:\manifold\to\mathbb{R}$ is a scalar observable.
For brevity, we refer to $\varphi_k$ and $\mu_k$, as the $k$-th Koopman
eigenfunction and eigenvalue respectively.
%
%%%
%%%  you say k-th, you mean countable -- Clancy ?
%%%
%
We also define $\lambda_{k}=\frac{\log(\mu_{k})}{\Delta t}$, which is the
discrete time approximation of the continuous time eigenvalue.
Accompanying the eigenfunctions and eigenvalues are the {\em Koopman modes},
which are vectors in $\mathbb{C}^N$ (or spatial profiles if the dynamical system
is a spatially extended one) that can be used to reconstruct (or predict) the full state when combined with the Koopman eigenfunctions and eigenvalues~\cite{Rowley2009,williams_submitted}.
The Koopman eigenvalues, eigenfunctions, and modes have been used in many
applications including the analysis of  fluid
flows~\cite{Schmid2010,Schmid2011,Schmid2012,Rowley2009,Seena2011},  power
systems~\cite{Susuki2013,Susuki2012,Susuki2011},  nonlinear
stability~\cite{Mauroy2013a}, and state space parameterization~\cite{Mauroy2013};  here we exploit their  properties for data fusion purposes.

For this application, the ability of the Koopman eigenfunctions to generate a \emph{parameterization of state space} is key.
In the simple example that follows we use the phase of an ``oscillatory'' eigenfunction, which has $|\mu_k|=1$, and the magnitude of a ``decaying'' eigenfunction, which has $|\mu_k| < 1$, 
as an intrinsic (quasi action-angle) coordinate system for the slow manifold
(i.e., the ``weakest'' stable manifold) of a limit cycle.
While there are many data driven methods for nonlinear manifold parameterization (see, e.g., Ref.~\cite{Lee2007,Nadler2005}), the benefit of this approach is that the resulting parameterization is, in principle, {\em invariant to invertible transformations of the underlying system} and, in that sense, are an intrinsic set of coordinates for the system.

Mathematically,  if it is possible to reconstruct the underlying system state
from one snapshot of observation data, then $\obsfuna$ formally has an inverse, 
which we denote as $\vec{T}$ (i.e., $\state=\vec{T}(\obsa)$).
When this is not the case naturally, one can sometimes construct an ``extended''
$\obsa$ where such a $\vec T$ does exist by combining measurements taken at the
current and a finite number of previous times~\cite{juang1994applied,
Chorin2000,Stengel2012}.
Then if $\varphi:\manifold\to\mathbb{C}$ is a Koopman eigenfunction,
%%%
%%%  talk with Clancy about this too
%%%
 $\tilde{\varphi}=\tilde{\alpha}\varphi\circ\vec{T}$, where
$\tilde{\varphi}:\Ma\to\mathbb{C}$, is formally an eigenfunction of the Koopman
operator  with the eigenvalue $\mu$ for one set of sensors (rather than for the
full state).
The constant $\tilde{\alpha}\in\mathbb{C}$ highlights that this
eigenfunction is only defined up to a constant.
%
%This is straightforward to show.
%
The evolution operator expressed in terms of $\obsa$ is $\vec{\tilde{F}}=\obsfuna\circ\vec{F}\circ\vec{T}$ and the action of the associated Koopman operator is $\tilde{\koopman}\psi=\psi\circ\vec{\tilde{F}}$.
Then, $(\tilde{\koopman}\tilde{\varphi})(\obsa)=(\tilde{\varphi}\circ\vec{\tilde{F}})(\obsa)=\tilde{\alpha}(\varphi\circ\vec{T}\circ\obsfuna\circ\vec{F}\circ\vec{T}\circ\obsfuna)(\vec{x})=\tilde{\alpha}(\varphi\circ\vec{F})(\vec{x})=\tilde{\alpha}(\koopman\varphi)(\vec{x})=\mu\tilde{\alpha}\varphi(\vec{x})=\mu\tilde{\varphi}(\obsa)$, where we have assumed that $\tilde{\varphi}$ is still an element of some space of scalar observables.
This same argument could be used to obtain a $\hat{\varphi}$ and $\hat{\alpha}$
for the measurements represented by $\obsb$.
Finally, we define the ratio $\alpha = \hat{\alpha}/\tilde{\alpha}$, whose role we
will explain shortly.

To approximate these quantities, we use Extended Dynamic Mode Decomposition
(EDMD)~\cite{williams_submitted},  which is a recently developed data-driven
method that approximates the Koopman  ``tuples'' (i.e., triplets of related
eigenvalues, eigenfunctions and modes).
%%%
%%% Yannis says, if you have space, the first time you say "modes" put a footnote that usually people associate
%%% modes with the eigenfunction "cf Fourier modes",, but YOU means something else     -- Clancy
%%%
%
The inputs to the EDMD procedure are sets of snapshot pairs,
$\{(\vec{x}_{m},\vec{y}_{m})\}_{m=1,\ldots,M}$, and a set of dictionary elements
that span a subspace of scalar observables, which we represent as the
vector valued function $\dictionary(\state)=[\psi_{1}(\state),\psi_{2}(\state),\ldots,\psi_{K}(\state)]^{T}$ where $\psi_{k}:\manifold\to\mathbb{R}$ and $\dictionary:\manifold\to\mathbb{R}^{K}$.
%%%
%%% do you want to say someghting, possibly in SI, about bases and dictionaries ?   Clancy
%%%
 %
 This procedure results in the matrix
\begin{equation}
 \edmdK\triangleq\mat{G}^{+}\mat{A}\in\mathbb{R}^{K\times K},
 \end{equation}
  which is a finite dimensional approximation of the Koopman operator, where $\mat{G}=\sum_{m=1}^{M}\dictionary(\vec{x}_{m})\dictionary(\vec{x}_{m})^{T}$, $\mat{A}=\sum_{m=1}^{M}\dictionary(\vec{x}_{m})\dictionary(\vec{y}_{m})^{T}$, and $+$ denotes the pseudo-inverse.
The $k$-th eigenvalue and eigenvector of $\edmdK$, which we denote as $\mu_k$ and $\vec\xi_k$, produce an approximation of the $k$-th eigenvalue and eigenfunction of the Koopman operator~\cite{williams_submitted}.
We denote the approximate eigenfunction as
$\varphi_k(\vec{x})=\dictionary^{T}(\vec{x})\vec{{\xi}}_{k} = \sum_{i=1}^K
\psi_i(\vec x)\vec{\xi}_k^{(i)}$ where $\vec\xi_k^{(i)}\in\mathbb{C}$ is the
$i$-th element of $\vec \xi_k$.

{\bf The numerical procedure}. 
Because $\state$ is, by assumption, unknown, we instead apply EDMD to the
measurement data.
The $\obsa$ data generates approximations of the set of $\tilde{\varphi}_{k}$ and $\tilde{\mu}_{k}$, and the $\obsb$ data generates approximations of the set of $\hat{\varphi}_{k}$ and $\hat{\mu}_{k}$.
To map between these separate sets of observations, note that
\begin{equation}
\tilde{\varphi}_{k}(\obsa_{m})=\frac{\hat{\alpha}_{k}}{\tilde{\alpha}_{k}}\hat{\varphi}_{k}(\obsb_{m})=\alpha_{k}\hat{\varphi}_{k}(\obsb_{m}),\label{eq:mapping}
\end{equation}
because $\varphi_{k}(\state_{m})=\tilde{\alpha}_{k}\tilde{\varphi}_{k}(\obsa_{m})=\hat{\alpha}_{k}\hat{\varphi}_{k}(\obsb_{m})$ when $\tilde{\varphi}_{k}$ is the eigenfunction that ``corresponds'' to $\hat{\varphi}_{k}$.
To determine which eigenfunctions correspond to one another, we require that $\tilde{\mu}_{k}\approx\hat{\mu}_{k}$.
This is also a sanity check that EDMD is indeed producing a reasonable
approximation of the Koopman operator; if no eigenvalues satisfy this
constraint,  then the approximation generated by EDMD is not accurate enough to be useful.
Finally, to determine the $\alpha_{k}$, we use the joint data set along with \eqref{eq:mapping} to solve for $\alpha_{k}$ in a least squares sense (and thus register the data).
Taken together, the steps above produce an approximation of $\tilde{\varphi}_k$ given a measurement of $\obsb$.
The final step in the procedure is to obtain a mapping from $\tilde{\varphi}_k$
to $\obsa$.
One conceptually appealing way is by expressing $\obsa$ as the sum of its Koopman modes and eigenfunctions~\cite{Rowley2009,Budivsic2012,williams_submitted}.
In this manuscript, we take a simpler approach and use interpolation routines for scattered data (see, e.g., Refs.~\cite{Lancaster1981,Levin1998}) to approximate the inverse mapping, $(\varphi_1(\obsa), \varphi_2(\obsa),\ldots) \mapsto \obsa$.
In particular, we use two-dimensional piecewise linear interpolation as
implemented by the \texttt{griddata} command in Scipy~\cite{jones2001scipy}.

Overall, the data merging procedure is as follows:
\begin{enumerate}
\itemsep0em
\item Approximate the eigenfunctions and eigenvalues with the data
$\{(\obsa_m,\obsaf_m)\}_{m=1,\ldots,\Ma}$ and
$\{(\obsb_{m},\obsbf_{m})\}_{m=1,\ldots,\Mb}$,  and determine the number of
eigenfunctions required to usefully parameterize state space.
\item Match the $\tilde{\varphi}_{k}$ with its corresponding
$\hat{\varphi}_{k}$ by finding the closest pair of eigenvalues, $\tilde\mu_k$ and $\hat\mu_k$.
\item Use the joint data set, $\left\{
(\obsa_{m},\obsb_{m})\right\}_{m=1}^{M_\text{joint}}$, to compute the
$\alpha_{k}$ for each eigenfunction.
\item Given a new measurement, $\obsb$, compute $\hat{\varphi}_{k}(\obsb)$ and
use \eqref{eq:mapping} to approximate $\tilde{\varphi}_{k}(\obsa)$.
\item Finally, approximate $\obsa$ from the
$\tilde{\varphi}_{k}(\obsa)$ using an interpolation routine.
\end{enumerate}
This method can be considered a hybrid of static and dynamic state estimation
techniques:
 dynamic information is required to construct the mapping to and from the set of
 intrinsic variables, but is not used beyond that.

{\bf An illustrative example.} To demonstrate the efficacy of this technique, we
apply it to the FitzHugh-Nagumo PDE in 1D, which is given by:
 \begin{subequations}
\label{eq:fhne}
\begin{align}
\partial_{t}v & =\partial_{xx}v+v-w-v^{3},\\
\partial_{t}w & =\delta\partial_{xx}w+\epsilon(v-c_{1}w-c_{0}),
\end{align}
\end{subequations}
where $v$ is the activation field, $w$ is the inhibition field, $c_{0}=-0.03$, $c_{1}=2.0$, $\delta=4.0$, $\epsilon=0.017$, for $x\in[0,20]$ with Neumann boundary conditions.
These parameter values are chosen so that \eqref{eq:fhne} has a spatio-temporal
limit cycle.
While the Koopman operator can be defined for~\eqref{eq:fhne},  the large dimension
of the state space (e.g., for a finite difference discretization) would make the necessary computations intractable.
Instead, the simpler task of constructing this mapping on a low-dimensional ``slow'' manifold, where the dynamics are effectively two-dimensional, is undertaken.

We start by collecting a large, representative data set, and performing
Principal Component Analysis (PCA)~\cite{Holmes1998} on it.
For the first set of measurements, we project $v$ and $w$ onto the leading {\em
three} principal  components of the complete data set, so
${\obsa_n=[a_{1}(t_n),a_{2}(t_n),a_{3}(t_n)]}$ where $a_{k}$ is the coefficient
of the $k$-th principal component evaluated at $t_n$.
Together these modes capture over 95\% of the energy in the data set~\cite{Holmes1998}, and serve as an accurate, low-dimensional representation of the full state.
The other data come from pointwise measurements taken at $x=10$  (i.e., $\obsb(t)=[v(10,t),w(10,t)]$).
Both sets of data are collected by generating 20 trajectories with a sampling interval of $\Delta t=2$, where each trajectory consists of  $10^{3}$ snapshot pairs.
Each trajectory is computed by perturbing the unstable fixed point associated with the limit cycle, evolving \eqref{eq:fhne} for $\Delta T = 1000$ (roughly ten periods of oscillation) to allow fast transients to decay, and then recording the evolution of each measurement.
The data sets are generated independently, \emph{so no snapshot in one data set
corresponds exactly to any of the snapshots in the other}, and then ``whitened''
so that each measurement set has unit variance.
The registration data set consists of a {\em single measurement} from the PCA data set and the corresponding pointwise values.

We approximate the space of observables with a finite dimensional dictionary,
whose elements we denote as $\tilde{\psi}_k$ or $\hat{\psi}_k$.
In this application, $\tilde{\psi}_{k}$ and $\hat{\psi}_{k}$ are the $k$-th
shape function of a moving least squares interpolant~\cite{Lancaster1981,Levin1998} 
with up to linear terms in each
variable~\cite{Li2002,Belytschko1994,Belytschko1996} and cubic spline weight
functions~\cite{Belytschko1996}.
Using a quad- or oct-tree to determine the node centers~\cite{Samet1990}
resulted in a set of 1622 basis functions for the PCA data and 1802 basis functions for the pointwise data.
Though non-trivial to implement, this set was chosen because it can exactly
reproduce constant and linear functions, but other choices of basis
functions are also viable~\cite{Wendland1999,Fasshauer1996,Karniadakis2013,Cockburn2000,Trefethen2000,Weideman2000,williams_submitted}.

\begin{figure*}[t]
\centering \includegraphics[width=0.9\textwidth]{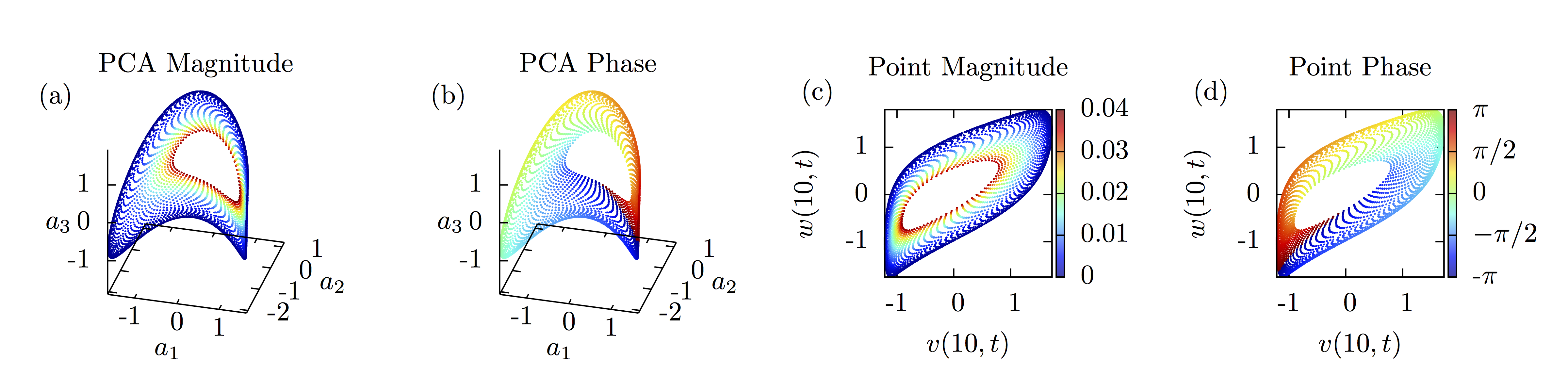}
\caption{The data driven parameterization generated using the Koopman eigenfunctions
with eigenvalues near $\lambda_{1}=-8\times10^{-4}$ and $\lambda_{2}=4.7i\times10^{-2}$
for the magnitude and phase figures respectively.
Plots (a) and (c)
show $\varphi_{1}$, the  value of the first eigenfunction,
and plots (b) and (d) show $\angle\varphi_{2}$ the ``phase'' of
the second eigenfunction. 
These pairs of eigenfunctions
can be use to parameterize the data; for the PCA case, the data effectively lie
on a nonlinear, 2-dimensional manifold embedded in $\mathbb{R}^{3}$,
and for the pointwise data, on a subset of $\mathbb{R}^{2}$.}
\label{fig:leading_eigenfunction_correspondence}
\end{figure*}

We applied the EDMD procedure to the dynamical data sets to obtain approximations of the Koopman eigenfunctions and eigenvalues.
 Figure~\ref{fig:leading_eigenfunction_correspondence} shows the PCA and point data colored by the magnitude of the Koopman eigenfunction with $\lambda_{1}\approx-8\times10^{-4}$ and by the phase of the Koopman eigenfunction with $\lambda_{2}\approx4.7i\times10^{-2}$.
In particular, $\tilde{\lambda}_{1}=-7.26\times 10^{-4}$,
$\hat{\lambda}_{1}=-8.57\times 10^{-4}$, $\tilde{\lambda}_{2}=0.0473i$, and $\hat{\lambda}_{2}=0.0473i$, where $\tilde{\lambda}_{k}$ is the $k$-th eigenvalue computed using the $\obsa$ measurements, and $\hat{\lambda}_{k}$ the eigenvalue obtained using $\obsb$ measurements.
Despite the differences in the nature of data, both relevant sets of Koopman
eigenfunctions generate (effectively) the same parameterization of the slow
manifold.

\begin{figure}[tb]
\centering \includegraphics[width=0.5\textwidth]{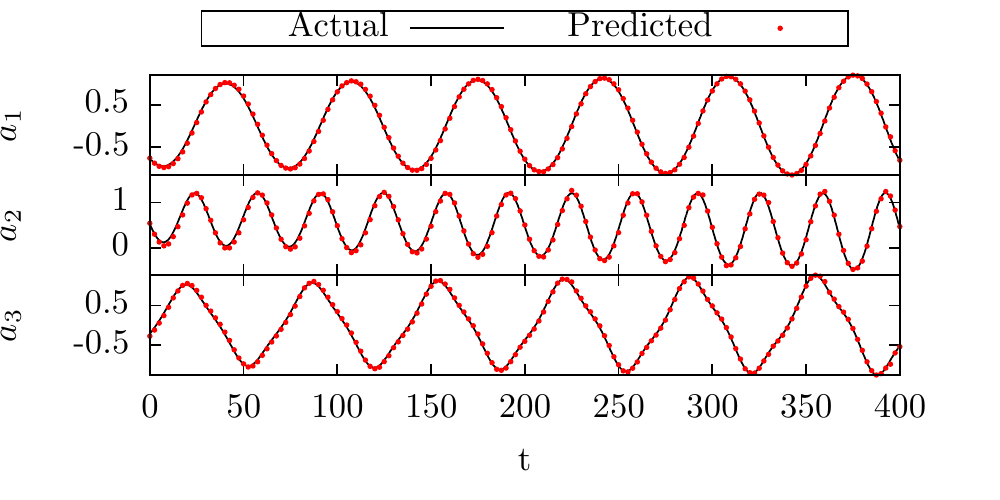}
\caption{The principal component coefficients reconstructed from the pointwise
data (red) for a \emph{new} trajectory that was not used in the Koopman eigenfunction
computation. 
The actual values are indicated by the black lines. 
There is good agreement between the predicted and true principal component
values, and our approach accurately captures the approach of the trajectory to
the limit cycle.
We reiterate that each principal component value was reconstructed
independently.
}
\label{fig:reconstruct}
\end{figure}

Using this pair of parameterizations, we reconstruct the PCA coefficients from pointwise data.
Figure~\ref{fig:reconstruct} demonstrates the quality of this reconstruction
with data from a randomly initialized trajectory that approaches the limit
cycle.
In the figure, the black lines denote the true coefficient of each principal component, and the red dots indicate the predicted value.
Note that this trajectory \emph{was not used to compute the Koopman
eigenfunctions}; furthermore, the fact that the data are a time series was not
used in making this prediction, and each measurement of $\obsb$ was considered individually.

Visually, the agreement between the predicted and actual states is good, but
there are  errors in the reconstruction.
Over the time window $t\in[0, 400]$, which is shown in
Fig.~\ref{fig:reconstruct}, the relative error in each of the principal
components is $e_{1}=0.0405$, $e_{2}=0.0655$, and $e_{3}=0.0496$, where
$e_{i}=\|a_{i}^{\text{(true)}}-a_{i}^{\text{(pred)}}\|/\|a_{i}^{(\text{true})}\|$,
where $\|\cdot\|=\left(\int_{0}^{400}(\cdot)^{2}\; dt\right)^{1/2}$.
In general, the accuracy of this approach is better for points nearer to the
limit cycle, and the relative errors over the window $t\in[0, 4000]$,
which (as shown in the supplement) contains points closer to the limit cycle,
are only $e_{1}=0.0140$, $e_{2}=0.0295$, and $e_{3}=0.0234$ respectively.

In this example, we have not yet discussed the dimensionality of the data, but
as stated previously, the number of measurements in each data set is critical
for this approach to be justifiable mathematically.
Our focus is on dynamics near the limit cycle, which in this problem are
effectively two dimensional.
However, the data lie on a two-dimensional {\em nonlinear manifold} that PCA,
which fits the data with a hyperplane, requires {\em three} principal
components to accurately represent.
Therefore, the identified mapping is from $\mathbb{R}^2$ to a two-dimensional
nonlinear manifold in $\mathbb{R}^3$, and not from $\mathbb{R}^2$ to
$\mathbb{R}^3$.

{\bf Conclusions.} We have presented a method for data
fusion or state reconstruction that is suitable for nonlinear systems and
exploits the existence of an intrinsic set of variables generated by the eigenfunctions
of the Koopman operator.
Rather than mapping directly between different sets of measurements, our method
focuses on generating an independent mapping to and from the intrinsic
variables for each heterogeneous set of measurements.
In principle, this can be accomplished as long as a mapping from each set of
measurements (or measurements and their time delays) to the system state exists,
and the benefit of this approach is that the majority of the required data can be
obtained independently, and only a single ``joint'' pair of data is needed.
The keys to this method are: (i) the invariance of the Koopman eigenvalues to
invertible transformations, (ii) the fact that the eigenfunctions parameterize
state space, and (iii) the ability of data-driven methods, such as EDMD, to
produce approximations of the eigenvalues and eigenfunctions that are ``accurate
enough'' to allow these properties to be exploited in practical settings.
%

%In the example presented above, we assumed the sampling interval, $\Delta t$,
%was the same for both sets of data,  but this could be relaxed with a small
% modification to the algorithm.
%
%With this modification, this approach could be used  to estimate time-resolved
%evolution of principal component  coefficients from an alternative, and easier
%to acquire, source of data similar to the study in Ref.~\cite{Tu2013}.
%
%Ultimately, while the Koopman operator is a powerful framework for the analysis
%of nonlinear  dynamical systems, the properties possessed by its eigenvalues
% and eigenfunctions also have practical applications; and although obtaining highly 
%accurate approximations of these quantities remains a challenge, data driven
%methods,  such as Extended Dynamic Mode Decomposition, are capable of producing 
%approximations that are accurate enough to be used for practical tasks such as 
%data merging and state reconstruction.
%

%

\textbf{Acknowledgements}  The authors gratefully acknowledge support from the
National Science Foundation (DMS-1204783) and the AFOSR.

\bibliographystyle{apsrev4-1}
\bibliography{koopman,sensor_placement}

\appendix
\section{Supplementary Material}

In this supplement, we present an expanded discussion of some key points in the
data fusion  procedure outlined in the manuscript.
In particular, we focus on the ability of the Koopman eigenfunctions to parameterize
state space, the accuracy of the mapping between eigenfunctions obtained with
different  sets of data, and the quality of the reconstruction for our
illustrative example  over a longer time horizon.
Each of these points has been touched upon in the manuscript, and our objective
is  simply to present some additional supporting evidence rather than introduce any new concepts.

\subsubsection*{Measurement Parameterization and Mapping }

\begin{figure*}
	\centering
	
	\includegraphics[width=0.95\textwidth]{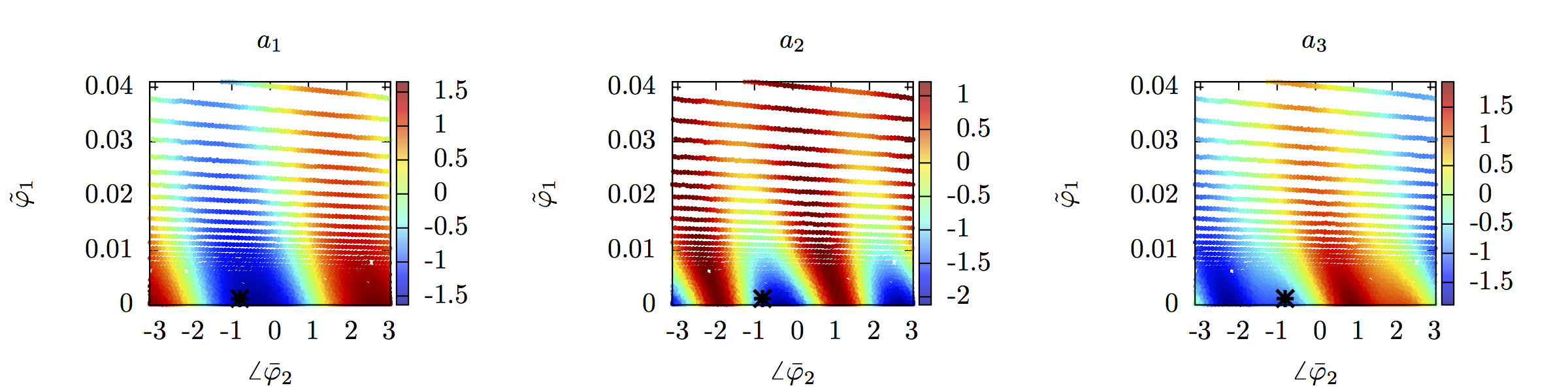}
	\caption{Pseudo-color plot of the first three principal component coefficients,
		$a_{1}$, $a_{2}$, and $a_{3}$, as a function of two of the Koopman
		eigenfunctions, $\tilde{\varphi}_{1}$ and $\angle\tilde{\varphi}_{2}$.
		Note that all three of these coefficients appear to be smooth functions of the Koopman eigenfunction values, and hence $\tilde{\varphi}_{1}$ and
		$\angle\tilde{\varphi}_{2}$ form an effective action-angle coordinate system for the nonlinear ``slow manifold" in principal component space.
		The black $x$ denotes the point in the joint data set that was used for
		registration purposes. }
	\label{fig:merging-pca}
\end{figure*}

In the manuscript, we claim that our set of two Koopman eigenfunctions
parameterizes both the PCA data and the point-wise measurement data,   and give
a visual example of this in Fig.~\ref{fig:leading_eigenfunction_correspondence}.
However, the final step in reconstructing the principal component coefficients
(or point-wise measurements)  is mapping from the intrinsic coordinates defined
by  the Koopman eigenfunctions to the original set of variables.
In principle, the Koopman modes would allow this to be done, and would
effectively act as coefficients in a truncated series expansion of the inverse map.
To simplify our code, we instead use interpolation methods appropriate for scattered data.
These methods could include moving least squares interpolants, but we make use
of the linear  interpolation routine implemented by the \texttt{griddata} command in Scipy.
This routine uses a Delaunay triangulation in conjunction with Barycentric
interpolation to create a piecewise linear approximation of the inverse
function.
Figure~\ref{fig:merging-pca} plots the coefficient of the $i$-th principal component as a function of the two Koopman eigenfunctions; 
note that the value of the principal component varies continuously and ``slowly,'' which is why the simple piecewise linear interpolant  employed here is 
sufficient to produce a good approximation of the mapping from the Koopman eigenfunction values to the principal component coefficients.
With far fewer data points or more complicated attractors, this mapping may be more difficult to satisfactorily approximate, and would require more sophisticated interpolation/extrapolation algorithms than what we used here.

\begin{figure}
	\centering
	
	\includegraphics[width=0.45\textwidth]{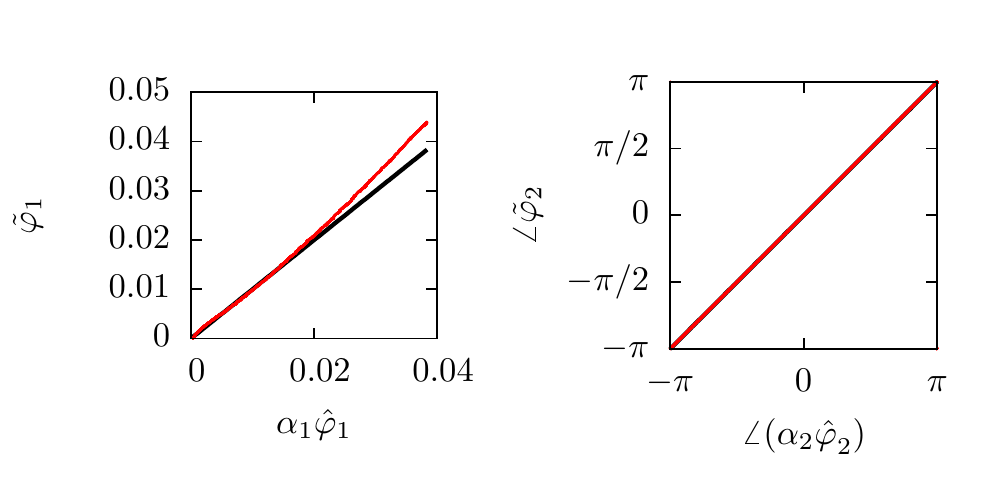}
	
	\caption{Plot of $\tilde{\varphi}_{1}$ as a function of $\alpha_{1}\hat{\varphi}_{1}$
		and $\angle\tilde{\varphi}_{2}$ as a function of $\angle(\alpha_{2}\hat{\varphi}_{2})$.
		The data points are indicated by the red dots; the black line is a
		guide for the eye. Ideally, the $\alpha$ are chosen to make these
		quantities equal. 
		While this is essentially the case for $\varphi_{2}$,
		the accuracy of the fit degrades as either $\tilde{\varphi}_{1}$ or
		$\hat{\varphi}_{1}$ become larger. 
		For this problem, our approximation of this
		eigenfunction is reasonably accurate when $|\tilde{\varphi}_{1}|<0.03$, and
		tends to be  higher for points nearer to the limit cycle (where $|\tilde\varphi_1|=0$).
		At larger values, the clear systemic error that appears in this plot will
		manifest itself as a systematic error in the reconstruction (or merging) of new
		data points. }
	\label{fig:one_to_one}
\end{figure}

A related issue involves the data used in the approximation. 
In Fig.~\ref{fig:merging-pca}, it is clear that the majority of the data lies near $\hat{\varphi}_{1}=0$.
There are two reasons for this.
First, it is easier to collect data near the limit cycle simply due to the dynamics of the underlying system.
Furthermore, the magnitude of the eigenfunction, $\varphi_{1}$, typically grows rapidly as one gets further from the limit cycle, and will even
have singularities at the unstable fixed point.
However, because the EDMD procedure uses data to approximate the eigenfunctions,
it is important to have data further away from the limit cycle despite the
(possible extra) effort in obtaining them.
It is for this reason that we use 20 independently initialized trajectories rather than a single, longer trajectory.

In the manuscript, we noted that the accuracy of our data fusion decreases the further we get from the limit cycle.
Figure~\ref{fig:one_to_one} quantifies this effect.
The figure shows the amplitude and phase of what ought to be the same eigenfunction (approximated through different sets of heterogeneous observations) as functions of one another.
Ideally, these two representations should coincide, and for $\angle\tilde{\varphi}_{2}$ and $\angle\hat{\varphi}_{2}$ they, in effect, do.
However, this is not the case for $\tilde{\varphi}_{1}$ and $\hat{\varphi}_{1}$,
which are shown in the left plot.
In theory, both of these eigenfunctions are zero on the periodic orbit, and their absolute value increases for points that are further away.
The figure, with the $\alpha$ we computed, shows good agreement between the
value of the eigenfunctions when their absolute values are small, but a growing
systematic error when they are large.
This difference is one of the causes of the errors we observe, and why the
accuracy of our approach  decreases as one gets further from the limit cycle.
There are, in general, open questions about the validity of EDMD computations
performed  on subsets of state space, and how they impact the accuracy
of  the numerically computed eigenfunctions.
For this problem, the predictions our method produces are accurate (to the eye)
when  $|\hat{\varphi}_{1}|<0.03$, but beyond that point, the predicted and
actual solutions will be quantitatively and visually different due to a
combination of sparse data and the ``partial domain'' issue.

\subsubsection*{Measurement Reconstruction}

\begin{figure*}
	\centering
	
	\includegraphics[width=0.9\paperwidth]{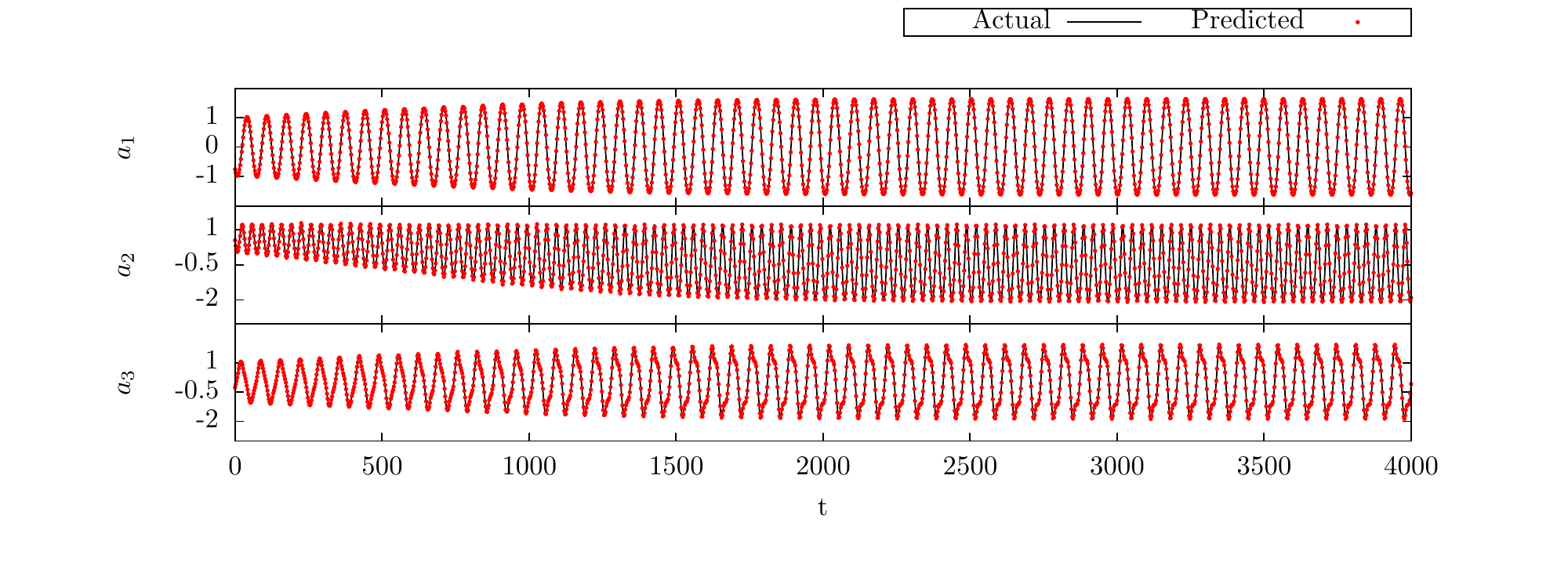}
	\caption{Plot of the predicted principal component coefficients obtained using
		the data merging procedure outlined in the paper. 
		This figure is a ``zoomed out'' version of Fig.~\ref{fig:reconstruct} and shows
		both the true and reconstructed values for $t\in[0, 4000]$, which highlights
		that this approach accurately recovers both the periodic oscillations and the
		``envelope'' that contains them. }
	\label{fig:long-prediction}
\end{figure*}

%%%% 
%%%% 800 is not clear unless you know zero 
%%%%
%%%% say the rules of the game 
%%%% and consider that I got distracted with the time - you really want to say
%%% you have better approximations close
%%%% 
%%%
%%% I submit you have a sentece to the effect of prediction
%%% because you can acttually predict 
%%% 
%%% do it to be inlcuded when it comes back 
%%% and we can then even say my favorite "precomputation" word !
%%%

In Fig.~\ref{fig:reconstruct}, we reconstructed the principal component values
from point-wise data for a new, randomly-initialized trajectory on the interval
$t\in[0, 400]$
Figure~\ref{fig:long-prediction} shows
the reconstructed state using the same trajectory, but over the time interval
$t\in[0,4000]$.
Overall, there is good agreement between the predicted and actual principal
component coefficients; in particular, our approach accurately recovers the
``envelope'' of the oscillation in all three principal components.
Although not obvious to the eye, the relative error decreases from $\sim$5\%
in the short window to $\sim$2\% over the longer time interval.
Recall that neither time nor the fact that this data are a trajectory are used
explicitly in the reconstruction; the increase in accuracy is solely due to the
fact that points at later times are near to the limit cycle.

\begin{figure*}
	\centering
	\includegraphics[width=0.8\textwidth]{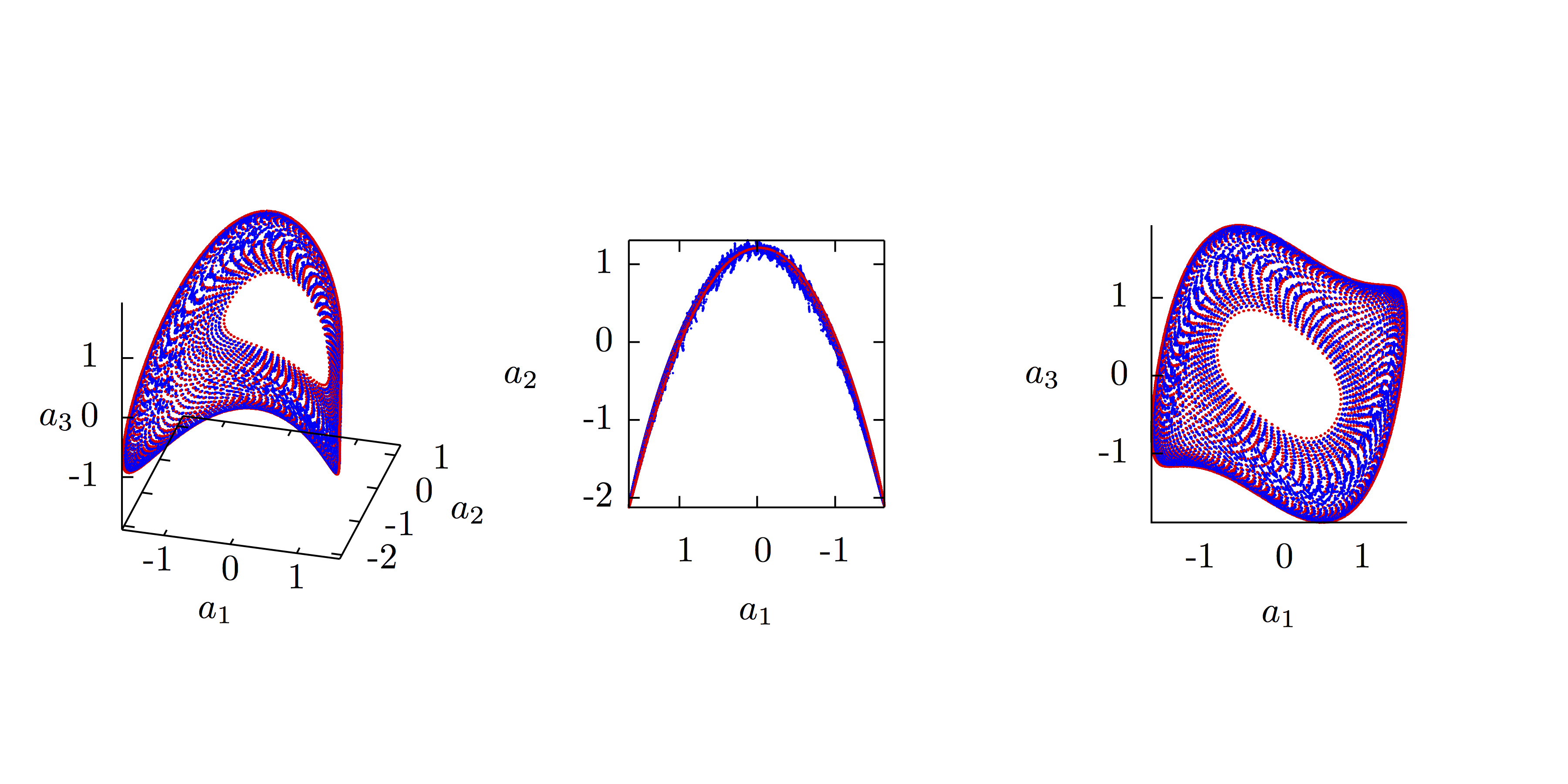}
	\caption{Three different views of the original principal component data  (shown
		in red), the reconstructed principal component values of the 20 trajectories that comprise
		the point-wise data set used in the manuscript (shown in blue), and the
		reconstruction with 20 new point-wise trajectories that were not used in the Koopman computation
		(also shown in blue).
		Though there is some error, which can  most easily be seen in the center plot,
		our approach effectively embeds these new points on the nonlinear,
		two-dimensional manifold that the original set of principal component
		coefficients is confined to.
		The region with the largest error occurs near the ``hole'' in the center of the manifold, and corresponds to points that are furthest from the limit cycle.
	}
	\label{fig:embded-new-manifold}
\end{figure*}

Finally, Fig.~\ref{fig:embded-new-manifold} shows the original set of PCA data
in blue, and 40 trajectories consisting  of the predicted values from the
point-wise reconstruction,  which are shown in red.
Twenty of the trajectories were the data used to construct the point-wise
approximation of the  Koopman eigenfunctions, and twenty of the trajectories are ``new''.
The three plots in the figure show the same data from three different views.
The point of this figure is to demonstrate that our predicted values lie on (or
near)  the nonlinear manifold on which the original PCA data live.
There are clearly some errors, which become particularly pronounced using the
view in the center  plot, that appear as ``noise'' about an otherwise parabolic shape.
As stated previously, these errors are most apparent near the ``hole'' in the
center of  the manifold, which corresponds to points far away from the limit
cycle where  our parameterizations are known to be inaccurate.
However, even these errors are small compared to the total change in the
coefficient  of the second principal component.
As a result, while this method is certainly not error free, it is able to
effectively map point-wise measurements to principal component values (and vice
versa) as long as the numerical approximations of the necessary Koopman
eigenfunctions are ``accurate enough'' to serve as a set of intrinsic variables.

\end{document}